\documentclass[a4paper, 12pt]{amsart}
\usepackage{amsmath, amsthm}
\usepackage{amssymb}
\usepackage{latexsym}
\usepackage{amsfonts}
\usepackage{graphicx}
\usepackage{graphics}
\usepackage{txfonts}
\newtheorem{thm}{Theorem}
\newtheorem{prop}[thm]{Proposition}
\newtheorem{lem}[thm]{Lemma}

\newtheorem{defi}{Definition}

\def\P{{\mathsf{P}}}

\def\F{{\mathcal{F}}}

\def\Z{ {\mathbf Z} }

\def\Ex{ \mathbf{E} }

\def\P{ \mathsf{P} }

\title{Local Time in Parisian Walkways}
\author{Jir\^o Akahori}
\address{
Department of Mathematics, Ritsumeikan University \\
1-1-1 Nojihigashi, Kusatsu, Shiga 525-8577, Japan}
\email{akahori@se.ritsumei.ac.jp}
\begin{document}
\begin{abstract}
In the present paper, Ito formula and Tanaka formula
for a special kind of symmetric random walk 
in the complex plane are studied.  
The random walk is called {\em Parisian walk}, and 
its {\em local time} is defined to be 
the number of exit from some regions.
\end{abstract}
\maketitle
\section{Parisian Walks}
Let $ \tau_1,...,\tau_n,... $ be an i.i.d. sequence with 
$ \P ( \tau =1 ) = \P( \tau = \zeta ) = \P ( \tau= \zeta^2 ) = 1/3 $,
where we denote $ \zeta = (-1 + \sqrt{-3})/2 $.
The filtration generated by $ \{\tau\} $ will be denoted by 
$ \mathbf{F} \equiv \{ \F_t \} $. 
\begin{defi}
An $ \mathbf{F}$-adapted complex valued process 
$ \{ Z_t \} $ is called {\bf Parisian} (walk) if
(i) it is a martingale starting from a point in 
$ \mathbb{Z}[\zeta] \equiv \{ a + b \zeta : a,b \in \mathbb{Z} \} $, 
and (ii) $ Z_{t+1}- Z_t \in \{ 1, \zeta, \zeta^2 \} $ 
for all $ t $.
\end{defi}
Thus, a Parisian walk is a random walk on $ \mathbb{Z}[ \zeta ] $. 
Note that there are a lot of Parisian walks as functions of $ \{ \tau \} $, 
but the law is unique up to the initial point.
\begin{figure}[h]
  \begin{center}
    \includegraphics[width=55mm,height=49mm]{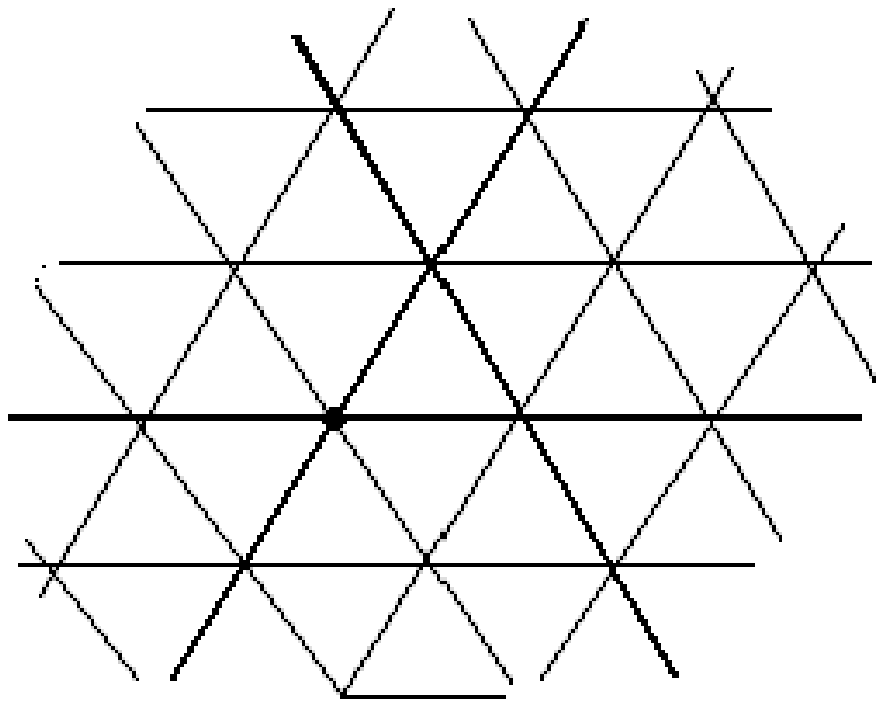}
\hspace{10mm}
    \includegraphics[width=57mm,height=49mm]{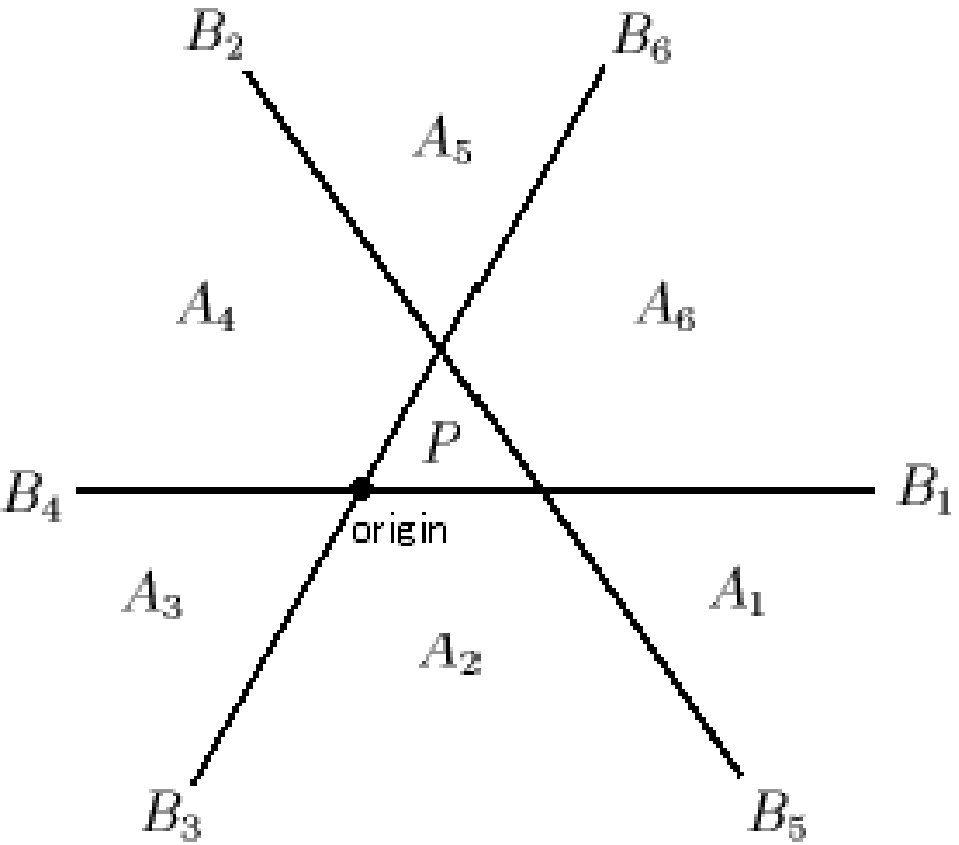}
  \end{center}
  \caption{Parisian walkways $ \mathbb{Z} [\zeta] $.}
  \label{fig:fig0.bmp}
\end{figure}

We classify the points in $ \mathbb{Z}[ \zeta ] $ as
\begin{equation*}
\begin{split}
& P = \{ 0,1,1+\zeta \}, \\
& A_1 = \{ k_1 + k_2 \zeta^2 \,|\, 1< k_2 +1< k_1 \},\, A_2 = \{ -k_1 \zeta^2 -k_2 \,|\, 0 \leq k_2 < k_1 \},\,\\
& A_3 = \{ -k_1 -k_2 \zeta \,|\, 0 < k_2 < k_1 \},\,
A_4 = \{ k_1 \zeta + k_2 \zeta^2 \,|\, 0 \leq k_2 < k_1 \},\, \\
& A_5 = \{ -k_1 \zeta^2 - k_2 \,|\, 1 < k_2 +1 < k_1 \},\,
A_6 = \{ k_1 + k_2 \zeta \,|\, 0 < k_2 < k_1 \},
\end{split}
\end{equation*}
\begin{equation*}
\begin{split}
& B_1 = \mathbb{Z}_{\geq 2},\, B_2 = \zeta \mathbb{Z}_{\geq 1}-\zeta^2,\, B_3 = \zeta^2 \mathbb{Z}_{\geq 1}, \\
&B_4 =  \mathbb{Z}_{\leq -1}, \, B_5 =  1+\zeta \mathbb{Z}_{\leq -1},
\, B_6 =  \zeta^2 \mathbb{Z}_{\leq -2}.
\end{split}
\end{equation*}
Further we define {\em closure}  of $ A_1, A_3 $ and $ A_5 $, respectively, as
\begin{equation*}
\begin{split}
& \bar{A}_1 := A_1 \cup B_1 \cup B_5 \cup \{ 1\} , \\
& \bar{A}_3 := A_3 \cup B_3 \cup B_4 \cup \{ 0 \} , \\
& \bar{A}_5 := A_5 \cup B_2 \cup B_6 \cup \{ 1 + \zeta \}.
\end{split}
\end{equation*}

In this paper. {\em local time} of a Parisian walk is defined to be the number of exit 
from $ \bar{A}_1 $, $ \bar{A}_3 $ and $ \bar{A}_5 $. More precisely, we set
\begin{equation*}
l^j_t := \sharp \{\, u \in Z_{ \geq 0} ,u < t \,|\, Z_u \in \bar{A}_j , \, Z_{u+1} \not\in \bar{A}_j \, \}, 
\end{equation*}
for $ j=1,3,5 $ and  {\bf local time} $ \{ L_t \} $ is defined by $ L_t = l^1_t + l^3_t + l^5_t $.

In addition, an $ \mathbf{F} $-martingale $ \{ X_t \} $ with $ X_0 \in \mathbb{Z} $ and
\begin{equation*}
X_t - X_{t-1} = 
\begin{cases}
1 & \text{with probability 1/3} \\
0 & \text{with probability 1/3} \\
-1 & \text{with probability 1/3}
\end{cases}
\end{equation*} 
will be called {\bf simple} (walk). 
Of course, the law of simple walks is unique up to the initial points.

In this paper, we present the following
\begin{thm}\label{0304152131p}
For a Parisian walk $ Z $, $ \Vert Z_t \Vert - L_t $ is a simple walk. 
Here, $ \Vert \cdot \Vert $ is the graph distance from $ P $ in $ \mathbb{Z} [\zeta] $,
namely it is the length of the shortest path from $ P $.
\end{thm}

\noindent {\bf Remark.}
Our results heavily rely on the arguments found in Fujita's paper \cite{Fuj} 
on 1-dimensional random walks. 

\section{An It\^o formula for Parisian walks}
We begin with the following lemma. 
\begin{lem}\label{0304301706p}
Let $ Z $ be a Parisian walk. 
Then the two dimensional process $ ( Z, \bar{Z} ) $ enjoys 
martingale representation property; 
every complex valued $ \mathbf{F} $-martingale is represented 
as a stochastic integral with respect to $ ( Z, \bar{Z} ) $.
\end{lem}
\begin{proof}
Denote $ \Delta Z_t := Z_t - Z_{t-1} $ for $ t \in \mathbb{Z}_{>0} $.
Fix $ t $ and set 
$$ \Delta Z_{S} := \prod_{s_i = \zeta } 
\Delta Z_i \prod_{s_i = \zeta^2} 
\overline{\Delta Z_i}$$ 
for $ S=(s_1,...,s_t) \in \{1, \zeta, \zeta^2 \}^t $. 
Then we have $ \Ex [\Delta Z_{S} \overline{ \Delta Z_{S'} }] = 1 $ 
if $ S= S' $ and $ =0 $ otherwise because of the martingale property
and of the fact that $ (\Delta Z_t)^2 = \overline{\Delta Z_t} $. 
Therefore $ \{ \Delta Z_S \,|\, S \in \{1, \zeta, \zeta^2 \}^t \} $ forms
an orthonormal basis (ONB) of $ L^2 (\F_t) $ 
since $ \sharp \{1, \zeta, \zeta^2 \}^t = \dim L^2 (\F_t)= 3^t $.

For an adapted $ \{ X_t \} $, expanding $ X_t - X_{t-1} $ with respect to
this ONB and denoting $ \Ex [(X_t-X_{t-1})\overline{\Delta Z_S} ] = x_S $, we have
\begin{equation}\label{0304112126p}
\begin{split}
& X_t - X_{t-1} \\
&= \sum_{s_t= \zeta}x_S \Delta Z_S 
+ \sum_{s_t=\zeta^2} x_S\Delta Z_S + \sum_{s_t=1} x_S \Delta Z_S \\
&= \left( \sum_{s_t= \zeta}x_S \Delta Z_{(s_1,...,s_{t-1})} \right) \Delta Z_t 
+ \left( \sum_{s_t=\zeta^2} x_S \Delta Z_{(s_1,...,s_{t-1})} \right) \Delta \bar{Z}_t
+ \sum_{s_t=1} x_S \Delta Z_{(s_1,...,s_{t-1})}.
\end{split}
\end{equation}
By summing up the above equation, we obtain the Doob decomposition of $ X $,
and this completes the proof.
\end{proof}

\noindent {\bf Remark.} 
The above lemma can be easily extended to
general unit root cases. The point here is that 
Parisian walk is the right discrete analogue of 
planer Brownian motion. 

\ 

Let $ \{ Z_t \} $ be a Parisian walk, 
and let $ f $ be a complex valued function on $ \mathbb{Z}[\zeta] $.
Then we have the following formula, which would correspond to an It\^o's formula
in $ \mathbf{F} $.
\begin{prop}
For $ t =0,1,2,...$, we have
\begin{equation}\label{0304112136p}
\begin{split}
f( Z_{t+1} ) - f (Z_t ) 
& = \frac{1}{3} (Z_{t+1} - Z_t)
\{ f(Z_t +1) + \zeta^2 f(Z_t + \zeta) + \zeta f(Z_t + \zeta^2) \} \\
&+ \frac{1}{3} ( \bar{Z}_{t+1} - \bar{Z}_t)
\{ f(Z_t +1) + \zeta f(Z_t + \zeta ) + \zeta^2 f(Z_t + \zeta^2 ) \} \\
&+ \frac{1}{3} 
\{ f(Z_t +1) + f(Z_t + \zeta ) + f(Z_t + \zeta^2 ) - 3 f(Z_t) \}.
\end{split}
\end{equation}
\end{prop}
\begin{proof}
As in the expression (\ref{0304112126p}), 
\begin{equation*}
f( Z_{t+1} ) - f (Z_t ) = \alpha \Delta Z_{t+1} + \beta \Delta \bar{Z}_{t+1} + \gamma
\end{equation*}
for some $ \F_t $-measurable $ \alpha, \beta $ and $ \gamma $. 
On the set of $ \Delta Z_{t+1} = 1 $, 
$ \Delta Z_{t+1} = \zeta $, and $ \Delta Z_{t+1} = \zeta^2 $ respectively, we have
\begin{equation}\label{0304112135p}
\begin{split}
f(Z_{t} +1 ) - f(Z_t) &= \alpha + \beta + \gamma, \\
f(Z_{t} + \zeta ) - f(Z_t) &= \alpha \zeta + \beta \zeta^2 + \gamma,  \\
\text{and}\,\,f(Z_{t} + \zeta^2 ) - f(Z_t) &= \alpha \zeta^2 + \beta \zeta + \gamma.
\end{split}
\end{equation}
Solving (\ref{0304112135p}) in terms of $ (\alpha, \beta,\gamma) $, 
we obtain (\ref{0304112136p}).
\end{proof}
\section{A Tanaka formula for Parisian walks}
For a Parisian walk $ Z $, set 
\begin{equation*}
\begin{split}
& \varphi_t := 
 1_{\{Z_t \in \bar{A}_1\} } 
+ \zeta 1_{ \{ Z_t \in \bar{A}_3 \}}+ \zeta^2 1_{ \{Z_t \in \bar{A}_5\} } \\
& \psi_t := 1_{\{Z_t \in A_6\}} + \zeta 1_{\{ Z_t \in A_4 \}} 
+ \zeta^2 1_{\{ Z_t \in A_2 \} }. 
\end{split}
\end{equation*}
Then, we have the following
\begin{prop}[A Tanaka formula for Parisian walks]
For $ t =0,1,2,...$, we have
\begin{equation}\label{0304152127p}
\Vert Z_{t+1} \Vert - \Vert Z_t \Vert 
= \frac{2}{3} \mathrm{Re} \left( (1-\zeta^2) (\varphi_t \Delta Z_t + \psi_t \Delta \bar{Z}_t) \right) 
+ L_{t+1} - L_t.
\end{equation}
\end{prop}
\begin{proof}
Applying It\^o formula (\ref{0304112136p}) to $ f (z) = \Vert z \Vert $, we have
\begin{equation*}
\begin{split}
\Vert Z_{t+1} \Vert - \Vert Z_t \Vert 
&= \frac{1}{3}
( \Vert Z_t + 1 \Vert + \zeta^2 \Vert Z_t + \zeta \Vert + \zeta \Vert Z_t + \zeta^2 \Vert )
\Delta Z_t \\
&+ \frac{1}{3}
( \Vert Z_t + 1 \Vert + \zeta \Vert Z_t + \zeta \Vert + \zeta^2 \Vert Z_t + \zeta^2 \Vert )
\Delta \bar{Z}_t \\
&+ \frac{1}{3}
( \Vert Z_t + 1 \Vert + \Vert Z_t + \zeta \Vert + \Vert Z_t + \zeta^2 \Vert 
- 3 \Vert Z_t \Vert).
\end{split}
\end{equation*}
Set 

\begin{equation*}
\begin{split}
 g_1 (z) &= \Vert z + 1 \Vert + \zeta^2 \Vert z 
+ \zeta \Vert + \zeta \Vert z + \zeta^2 \Vert \\
g_2 (z) &= \Vert z + 1 \Vert + \zeta \Vert z 
+ \zeta \Vert + \zeta^2 \Vert z + \zeta^2 \Vert \\
g_3(z) &=  \Vert z + 1 \Vert + \Vert z 
+ \zeta \Vert + \Vert z + \zeta^2 \Vert - 3 \Vert z \Vert. 
\end{split}
\end{equation*}
Then, we have the following tables.
\begin{table}[h]
 \caption{}
 {\tiny
 \begin{center}

  \begin{tabular}{|c|c|c|c|c|c|c|c|c|c|c|c|c|}
    \hline
       & $ A_1 $ & $ A_2 $ & $ A_3 $ & $ A_4 $ & $ A_5$ & $ A_6 $ &
 $ B_1 $ & $ B_2 $ & $ B_3 $ & $ B_4 $ & $ B_5 $ & $ B_6 $  \\
    \hline
    $ g_1 $ & $ 1-\zeta^2 $ & $ \zeta - \zeta^2 $ & $ \zeta-1 $ & 
$ \zeta^2 -1 $ & $\zeta^2 -\zeta $ & $ 1-\zeta $ & 
$ 1 $ & $ \zeta^2 $ & $ \zeta $ & $  -2 $ & $ -2\zeta^2 $ & $ -2\zeta $ \\
    \hline
    $ g_2 $ & $ 1-\zeta $ & $ \zeta^2 -\zeta $ & $ \zeta^2 -1 $ &
 $ \zeta -1 $ & $ \zeta -\zeta^2 $ & $1- \zeta^2$ & 
$1$ & $ \zeta $ & $\zeta^2$ & $ -2 $ & $-2\zeta$ & $  -2\zeta^2 $ \\
    \hline
    $ g_3 $ & $ 0 $ & $ 0 $ & $0$ & $0$ & $0$ & $0$ 
& $ 1 $ & $1$ & $1$ & $1$ & $1$ & $1$ \\
   \hline
  \end{tabular}
 \end{center}
 }
\end{table}
\begin{table}[h]
\caption{}
\begin{center}
\begin{tabular}{|c|c|c|c|}
\hline
 & $ 0 $ & $ 1 $ & $ 1 + \zeta $ \\
\hline
$ g_1 $ &$ -1 $ & $ -\zeta^2 $ & $ -\zeta $ \\
\hline
$ g_2 $ & $-1$ & $ -\zeta $ & $ -\zeta^2 $ \\
\hline
$ g_3 $ & $ 2 $ & $ 2 $ & $ 2 $ \\
\hline 
\end{tabular}
\end{center}
\end{table}

Using these tables, we obtain, for example,
\begin{equation*}
\Vert Z_{t+1} \Vert - \Vert Z_t \Vert =  \frac{1}{3} ( \Delta Z_t + \Delta \bar{Z}_t + 1)
= \frac{1}{3} \left\{ (1-\zeta^2) \Delta Z_t + ( 1-\zeta) \Delta \bar{Z}_t \right\}  + 1_{ \{  \Delta Z_t =\zeta \}}
\end{equation*}
on the set $ \{ Z_t \in B_1 \} $, and similar expressions on $ B_2,...,B_6 $ and $ \{0 \}, \{ 1\}, \{ 1+ \zeta \} $ 
lead to the formula (\ref{0304152127p}).
\end{proof}

Now Theorem \ref{0304152131p} follows easily from this version of Tanaka formula, 
by observing that $ \sum ( \varphi \Delta Z + \psi \Delta \bar{Z} ) $ is a Parisian walk, and 
real part of $ 2 (1-\zeta^2) Z'/3 $ is a simple walk whenever $ Z' $ is Parisian.

\ 

\noindent {\bf Remark.}
As we have seen in the proof of Lemma \ref{0304301706p}, $ Z $ and $ \bar{Z} $ are
mutually orthogonal martingale so that, after taking appropriate scaling limit, $ Z $ 
converges to a so-called {\em planer} Brownian motion.
Then, we can prove the standard It\^o's formula from our It\^o's formula (\ref{0304112136p}), noting that
(\ref{0304112136p}) can be rewritten as follows if $ f $ is real valued.
\begin{equation*}
\Delta f = \mathrm{Re} (Df) \mathrm{Re} (\Delta Z) + \mathrm{Im} (Df) \mathrm{Im} (\Delta Z)
+ L f,
\end{equation*}
where we have set $ Lf = \sum_{j=0,1,2} (f(z+\zeta^j) - f(z))/3 $ so that $ L $ corresponds to
the {\em Discrete Laplacian} operator in the context of random walks on graphs,
and $ Df = \sum_{j=0,1,2} \zeta^j f(z+\zeta^j)/3 $, which corresponds to
a {\em differential} operator on $ \Z [\zeta] $.

Moreover, a continuous time version of our theorem is given in the following way.

\ 

\noindent {\sl Let $ (B^1 B^2) $ be a 2-dimensional Brownian motion. Then
there exists an increasing process $ \{ L \} $ which increases only on (continuous time limit of) 
the lines 
$B_1,...,B_6 $ such that 
\begin{equation*}
\max ( |B^2_t|, | B_t^1- B^2/\sqrt{3}|, |B^1_t + B^2/\sqrt{3}| ) - L_t \quad \text{is a standard Brownian motion}.
\end{equation*}
}

\noindent This fact is verified by using the standard Tanaka formula again and again. 

\ 

\noindent {\bf Remark.} 
We can extend the discrete It\^o formula 
to more general cases from a perspective of Fourier expansion.
See \cite{A}.


\begin{thebibliography}{abcd}
\bibitem[A]{A}
Akahori, J, ``Discrete Ito Formulas and Their Applications to Stochastic Numerics", RIMS kokyuroku 1462 202-210 (2006). arXiv math.PR/0603341

\bibitem[F]{Fuj}
Fujita, T. ``A random walk analogue of L\'evy's theorem.'' preprint, 2003.
\end{thebibliography}
\end{document}